# ASYMPTOTIC NORMALITY OF THE $L_K$-ERROR OF THE GRENANDER ESTIMATOR

By Vladimir N. Kulikov and Hendrik P. Lopuhaä

*Eurandom and Delft University of Technology*

We investigate the limit behavior of the $L_k$-distance between a decreasing density $f$ and its nonparametric maximum likelihood estimator $\hat{f}_n$ for $k \geq 1$. Due to the inconsistency of $\hat{f}_n$ at zero, the case $k = 2.5$ turns out to be a kind of transition point. We extend asymptotic normality of the $L_1$-distance to the $L_k$-distance for $1 \leq k < 2.5$, and obtain the analogous limiting result for a modification of the $L_k$-distance for $k \geq 2.5$. Since the $L_1$-distance is the area between $f$ and $\hat{f}_n$, which is also the area between the inverse $g$ of $f$ and the more tractable inverse $U_n$ of $\hat{f}_n$, the problem can be reduced immediately to deriving asymptotic normality of the $L_1$-distance between $U_n$ and $g$. Although we lose this easy correspondence for $k > 1$, we show that the $L_k$-distance between $f$ and $\hat{f}_n$ is asymptotically equivalent to the $L_k$-distance between $U_n$ and $g$.

**1. Introduction.** Let $f$ be a nonincreasing density with compact support. Without loss of generality, assume this to be the interval $[0, 1]$. The nonparametric maximum likelihood estimator $\hat{f}_n$ of $f$ was discovered by Grenander [2]. It is defined as the left derivative of the least concave majorant (LCM) of the empirical distribution function $F_n$ constructed from a sample $X_1, \ldots, X_n$ from $f$. Prakasa Rao [11] obtained the earliest result on the asymptotic pointwise behavior of the Grenander estimator. One immediately striking feature of this result is that the rate of convergence is of the same order as the rate of convergence of histogram estimators, and that the asymptotic distribution is not normal. It took much longer to develop distributional theory for global measures of performance for this estimator. The first distributional result for a global measure of deviation was the convergence to a normal distribution of the $L_1$-error mentioned in [3] (see [4] for a rigorous proof ). A similar result in the regression setting has been obtained by Durot [1].

In this paper we extend the result for the $L_1$-error to the $L_k$-error, for $k \geq 1$. We will follow the same approach as in [4], which, instead of comparing $\hat{f}_n$ to $f$, compared both inverses. The corresponding $L_1$-errors are the same, since they represent the area between the graphs of $\hat{f}_n$ and $f$ and the area between the graphs of the inverses. Clearly, for $k > 1$ we no longer have such an easy correspondence between the two $L_k$-errors. Nevertheless, we will show that the $L_k$-error between $\hat{f}_n$ and $f$ can still be approximated by a scaled version of the $L_k$-error between the two inverses and that this scaled version is asymptotically normal.

The main reason to do a preliminary inversion step is that we use results from [4] on the inverse process. But apart from this, we believe that working with $\hat{f}_n$ directly will not make life easier. For $a \in [f(1), f(0)]$, the (left continuous) inverse of $\hat{f}_n$ is $U_n(a) = \sup\{x \in [0,1] : \hat{f}_n(x) \geq a\}$. Since $\hat{f}_n(x)$ is the left continuous slope of the LCM of $F_n$ at the point $x$, a simple picture shows that it has the more useful representation

$$(1.1) \qquad U_n(a) = \underset{x \in [0,1]}{\arg\max} \{F_n(x) - ax\}.$$

Here the arg max function is the supremum of the times at which the maximum is attained. Since $U_n(a)$ can be seen as the $x$-coordinate of the point that is touched first when dropping a line with slope $a$ on $F_n$, with probability one $\hat{f}_n(x) \leq a$ if and only if $U_n(a) \leq x$. Asymptotic normality of the $L_k$-error relies on embedding the process $F_n(x) - ax$ into a Brownian motion with drift. The fact that the difference between $F_n(x) - ax$ and the limit process is small directly implies that the difference of the locations of their maxima is small. However, it does not necessarily imply that the difference of the slopes of the LCMs of both processes is small. Similarly, convergence in distribution of suitably scaled finite-dimensional projections of $U_n$ follows immediately from distributional convergence of $F_n$, after suitable scaling, and an arg max type of continuous mapping theorem (see, e.g., [6]). When working with $\hat{f}_n$ directly, similar to Lemma 4.1 in [11], one needs to bound the probability that the LCM of a Gaussian approximation of $F_n$ on $[0,1]$ differs from the one restricted to a shrinking interval, which is somewhat technical and tedious.

Another important difference between the case $k > 1$ and the case $k = 1$ is the fact that, for large $k$, the inconsistency of $\hat{f}_n$ at zero, as shown by Woodroofe and Sun [13], starts to dominate the behavior of the $L_k$-error. By using results from [9] on the behavior of $\hat{f}_n$ near the boundaries of the support of $f$, we will show that, for $1 \leq k < 2.5$, the $L_k$-error between $\hat{f}_n$ and $f$ is asymptotically normal. This result can be formulated as follows. Define, for $c \in \mathbb{R}$,

$$V(c) = \arg\max_{t \in \mathbb{R}} \{W(t) - (t-c)^2\}, \tag{1.2}$$

$$\xi(c) = V(c) - c, \tag{1.3}$$

where $\{W(t): -\infty < t < \infty\}$ denotes standard two-sided Brownian motion on $\mathbb{R}$ originating from zero [i.e., $W(0) = 0$].

THEOREM 1.1 (Main theorem). *Let $f$ be a decreasing density on $[0, 1]$ satisfying:*

(A1) $0 < f(1) \leq f(y) \leq f(x) \leq f(0) < \infty$, *for $0 \leq x \leq y \leq 1$;*
(A2) *$f$ is twice continuously differentiable;*
(A3) $\inf_{x \in (0,1)} |f'(x)| > 0$.

*Then for $1 \leq k < 2.5$, with $\mu_k = \{E|V(0)|^k \int_0^1 (4f(x)|f'(x)|)^{k/3}\, dx\}^{1/k}$, the random variable*

$$n^{1/6}\left\{n^{1/3}\left(\int_0^1 |\hat{f}_n(x) - f(x)|^k\, dx\right)^{1/k} - \mu_k\right\}$$

*converges in distribution to a normal random variable with zero mean and variance*

$$\frac{\int_0^1 f(x)^{(2k+1)/3}|f'(x)|^{(2k-2)/3}\, dx}{k^2(E|V(0)|^k \int_0^1 (f(x)|f'(x)|)^{k/3}\, dx)^{(2k-2)/k}} \cdot 8\int_0^\infty \mathrm{cov}(|\xi(0)|^k, |\xi(c)|^k)\, dc.$$

Note that the theorem holds under the same conditions as in [4]. For $k \geq 2.5$, Theorem 1.1 is no longer true. However, the results from [9] enable us to show that an analogous limiting result still holds for a modification of the $L_k$-error.

In Section 2 we introduce a Brownian approximation of $U_n$ and derive asymptotic normality of a scaled version of the $L_k$-distance between $U_n$ and the inverse $g$ of $f$. In Section 3 we show that on segments $[s, t]$ where the graph of $\hat{f}_n$ does not cross the graph of $f$, the difference

$$\left|\int_s^t |\hat{f}_n(x) - f(x)|^k\, dx - \int_{f(t)}^{f(s)} \frac{|U_n(a) - g(a)|^k}{|g'(a)|^{k-1}}\, da\right|$$

is of negligible order. Together with the behavior near the boundaries of the support of $f$, for $1 \leq k < 2.5$ we establish asymptotic normality of the $L_k$-distance between $\hat{f}_n$ and $f$ in Section 4. In Section 5 we investigate the case $k > 2.5$ and prove a result analogous to Theorem 1.1 for a modified $L_k$-error.

REMARK 1.1.   With almost no additional effort, one can establish asymptotic normality of a weighted $L_k$-error $n^{k/3}\int_0^1 |\hat{f}_n(t) - f(t)|^k w(t)\,dt$, where $w$ is continuously differentiable on $[0,1]$. This may be of interest when one wants to use weights proportional to negative powers of the limiting standard deviation $(\frac{1}{2}f(t)|f'(t)|)^{1/3}$ of $\hat{f}_n(t)$. Moreover, when $w$ is estimated at a sufficiently fast rate, one may also replace $w$ by its estimate in the above integral. Similar results are in [8] for a weighted $L_k$-error.

**2. Brownian approximation.**   In this section we will derive asymptotic normality of the $L_k$-error of the inverse process of the Grenander estimator. For this we follow the same line of reasoning as in Sections 3 and 4 in [4]. Therefore, we only mention the main steps and transfer all proofs to the Appendix.

Let $E_n$ denote the empirical process $\sqrt{n}(F_n - F)$. For $n \geq 1$, let $B_n$ be versions of the Brownian bridge constructed on the same probability space as the uniform empirical process $E_n \circ F^{-1}$ via the Hungarian embedding, and define versions $W_n$ of Brownian motion by

$$(2.1) \qquad W_n(t) = B_n(t) + \xi_n t, \qquad t \in [0,1],$$

where $\xi_n$ is a standard normal random variable, independent of $B_n$. For fixed $a \in (f(1), f(0))$ and $J = E, B, W$, define

$$(2.2) \qquad V_n^J(a) = \arg\max_t \{X_n^J(a,t) + n^{2/3}[F(g(a) + n^{-1/3}t) - F(g(a)) - n^{-1/3}at]\},$$

where

$$(2.3) \qquad \begin{aligned} X_n^E(a,t) &= n^{1/6}\{E_n(g(a) + n^{-1/3}t) - E_n(g(a))\}, \\ X_n^B(a,t) &= n^{1/6}\{B_n(F(g(a) + n^{-1/3}t)) - B_n(F(g(a)))\}, \\ X_n^W(a,t) &= n^{1/6}\{W_n(F(g(a) + n^{-1/3}t)) - W_n(F(g(a)))\}. \end{aligned}$$

One can easily check that $V_n^E(a) = n^{1/3}\{U_n(a) - g(a)\}$. A graphical interpretation and basic properties of $V_n^J$ are provided in [4]. For $n$ tending to infinity, properly scaled versions of $V_n^J$ will behave as $\xi(c)$ defined in (1.3).

As a first step, we prove asymptotic normality for a Brownian version of the $L_k$-distance between $U_n$ and $g$. This is an extension of Theorem 4.1 in [4].

THEOREM 2.1.   *Let $V_n^W$ be defined as in* (2.2) *and $\xi$ by* (1.3). *Then for $k \geq 1$,*

$$n^{1/6}\int_{f(1)}^{f(0)} \frac{|V_n^W(a)|^k - E|V_n^W(a)|^k}{|g'(a)|^{k-1}}\,da$$

*converges in distribution to a normal random variable with zero mean and variance*

$$\sigma^2 = 2\int_0^1 (4f(x))^{(2k+1)/3} |f'(x)|^{(2k-2)/3}\,dx \int_0^\infty \operatorname{cov}(|\xi(0)|^k, |\xi(c)|^k)\,dc.$$

The next lemma shows that the limiting expectation in Theorem 2.1 is equal to

$$\mu_k = \left\{ E|V(0)|^k \int_0^1 (4f(x)|f'(x)|)^{k/3}\,dx \right\}^{1/k}. \tag{2.4}$$

LEMMA 2.1. *Let $V_n^W$ be defined by (2.2) and let $\mu_k$ be defined by (2.4). Then for $k \geq 1$,*

$$\lim_{n \to \infty} n^{1/6} \left\{ \int_{f(1)}^{f(0)} \frac{E|V_n^W(a)|^k}{|g'(a)|^{k-1}}\,da - \mu_k^k \right\} = 0.$$

The next step is to transfer the result of Theorem 2.1 to the $L_k$-error of $V_n^E$. This can be done by means of the following lemma.

LEMMA 2.2. *For $J = E, B, W$, let $V_n^J$ be defined as in (2.2). Then for $k \geq 1$, we have*

$$n^{1/6} \int_{f(1)}^{f(0)} (|V_n^B(a)|^k - |V_n^W(a)|^k)\,da = o_p(1)$$

*and*

$$\int_{f(1)}^{f(0)} \big||V_n^E(a)|^k - |V_n^B(a)|^k\big|\,da = \mathcal{O}_p(n^{-1/3}(\log n)^{k+2}).$$

From Theorem 2.1 and Lemmas 2.1 and 2.2, we immediately have the following corollary.

COROLLARY 2.1. *Let $U_n$ be defined by (1.1) and let $\mu_k$ be defined by (2.4). Then for $k \geq 1$,*

$$n^{1/6}\left( n^{k/3} \int_{f(1)}^{f(0)} \frac{|U_n(a) - g(a)|^k}{|g'(a)|^{k-1}}\,da - \mu_k^k \right)$$

*converges in distribution to a normal random variable with zero mean and variance $\sigma^2$ defined in Theorem 2.1.*

**3. Relating both $L_k$-errors.** When $k = 1$, the $L_k$-error has an easy interpretation as the area between two graphs. In that case $\int |U_n(a) - g(a)|\, da$ is the same as $\int |\hat{f}_n(x) - f(x)|\, dx$, up to some boundary effects. This is precisely Corollary 2.1 in [4]. In this section we show that a similar approximation holds for $\int_s^t |\hat{f}_n(x) - f(x)|^k\, dx$ on segments $[s, t]$ where the graphs of $\hat{f}_n$ and $f$ do not intersect. In order to avoid boundary problems, we will apply this approximation in subsequent sections to a suitable cut-off version $\tilde{f}_n$ of $\hat{f}_n$.

LEMMA 3.1. *Let $\tilde{f}_n$ be a piecewise constant left-continuous nonincreasing function on $[0, 1]$ with a finite number of jumps. Suppose that $f(1) \leq \tilde{f}_n \leq f(0)$, and define its inverse function by*

$$\tilde{U}_n(a) = \sup\{x \in [0, 1] : \tilde{f}_n(x) \geq a\},$$

*for $a \in [f(1), f(0)]$. Suppose that $[s, t] \subseteq [0, 1]$, such that one of the following situations applies:*

1. $\tilde{f}_n(x) \geq f(x)$, *for* $x \in (s, t)$, *such that* $\tilde{f}_n(s) = f(s)$ *and* $\tilde{f}_n(t+) \leq f(t)$.
2. $\tilde{f}_n(x) \leq f(x)$, *for* $x \in (s, t)$, *such that* $\tilde{f}_n(t) = f(t)$ *and* $\tilde{f}_n(s) \geq f(s)$.

*If*

(3.1) $$\sup_{x \in [s,t]} |\tilde{f}_n(x) - f(x)| < \frac{(\inf_{x \in [0,1]} |f'(x)|)^2}{2 \sup_{x \in [0,1]} |f''(x)|},$$

*then for $k \geq 1$,*

$$\left| \int_s^t |\tilde{f}_n(x) - f(x)|^k\, dx - \int_{f(t)}^{f(s)} \frac{|\tilde{U}_n(a) - g(a)|^k}{|g'(a)|^{k-1}}\, da \right|$$

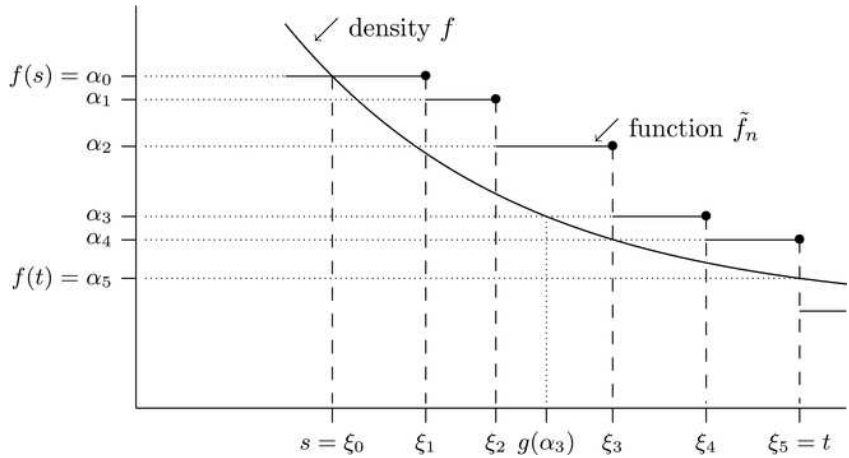

FIG. 1. *Segment $[s, t]$ where $\tilde{f}_n \geq f$.*

$$\leq C \int_{f(t)}^{f(s)} \frac{|\tilde{U}_n(a) - g(a)|^{k+1}}{|g'(a)|^k} \, da,$$

where $C > 0$ depends only on $f$ and $k$.

PROOF. Let us first consider case 1. Let $\tilde{f}_n$ have $m$ points of jump on $(s, t)$. Denote them in increasing order by $\xi_1 < \cdots < \xi_m$, and write $s = \xi_0$ and $\xi_{m+1} = t$. Denote by $\alpha_1 > \cdots > \alpha_m$ the points of jump of $\tilde{U}_n$ on the interval $(f(t), f(s))$ in decreasing order, and write $f(s) = \alpha_0$ and $\alpha_{m+1} = f(t)$ (see Figure 1). We then have

$$\int_s^t |\tilde{f}_n(x) - f(x)|^k \, dx = \sum_{i=0}^m \int_{\xi_i}^{\xi_{i+1}} |\tilde{f}_n(\xi_{i+1}) - f(x)|^k \, dx.$$

Apply a Taylor expansion to $f$ in the point $g(\alpha_i)$ for each term, and note that $\tilde{f}_n(\xi_{i+1}) = \alpha_i$. Then, if we abbreviate $g_i = g(\alpha_i)$ for $i = 0, 1, \ldots, m$, we can write the right-hand side as

$$\sum_{i=0}^m \int_{\xi_i}^{\xi_{i+1}} |f'(g_i)|^k (x - g_i)^k \left| 1 + \frac{f''(\theta_i)}{2f'(g_i)} (x - g_i) \right|^k dx$$

for some $\theta_i$ between $x$ and $g_i$, also using the fact that $g_i < \xi_i < x \leq \xi_{i+1}$. Due to condition (3.1) and the fact that $\tilde{f}_n(\xi_{i+1}) = \tilde{f}_n(x)$, for $x \in (\xi_i, \xi_{i+1}]$, we have that

(3.2)
$$\left| \frac{f''(\theta_i)}{f'(g_i)} (x - g_i) \right| \leq \frac{\sup |f''|}{\inf |f'|} \frac{|f(x) - f(g_i)|}{\inf |f'|}$$
$$\leq \frac{\sup |f''|}{(\inf |f'|)^2} |f(x) - \tilde{f}_n(x)| \leq \frac{1}{2}.$$

Hence, for $x \in (\xi_i, \xi_{i+1}]$

$$\left| \left| 1 + \frac{f''(\theta_i)(x - g_i)}{2f'(g_i)} \right|^k - 1 \right| \leq k \left( \frac{5}{4} \right)^{k-1} \frac{\sup |f''|}{2 \inf |f'|} (x - g_i).$$

Therefore, we obtain the inequality

$$\left| \int_s^t |\tilde{f}_n(x) - f(x)|^k \, dx - \sum_{i=0}^m \int_{\xi_i}^{\xi_{i+1}} |f'(g_i)|^k (x - g_i)^k \, dx \right|$$
$$\leq C_1 \sum_{i=0}^m \int_{\xi_i}^{\xi_{i+1}} (x - g_i)^{k+1} \, dx,$$

where $C_1 = k \sup |f'|^k (5/4)^{k-1} (\sup |f''|) / (2 \inf |f'|)$, or after integration,

$$\left| \int_s^t |\tilde{f}_n(x) - f(x)|^k \, dx \right.$$

$$(3.3) \quad \left. - \frac{1}{k+1} \sum_{i=0}^{m} |f'(g_i)|^k \{(\xi_{i+1} - g_i)^{k+1} - (\xi_i - g_i)^{k+1}\} \right|$$

$$\leq \frac{C_1}{k+2} \sum_{i=0}^{m} \{(\xi_{i+1} - g_i)^{k+2} - (\xi_i - g_i)^{k+2}\}.$$

Next consider the corresponding integral for the inverse $\tilde{U}_n$. Since $g_i < x < g_{i+1} < \xi_{i+1}$, we can write

$$\int_{f(t)}^{f(s)} \frac{|\tilde{U}_n(a) - g(a)|^k}{|g'(a)|^{k-1}}\, da = \sum_{i=0}^{m} \int_{\alpha_{i+1}}^{\alpha_i} \frac{|\xi_{i+1} - g(a)|^k}{|g'(a)|^{k-1}}\, da$$

$$= \sum_{i=0}^{m} \int_{g_i}^{g_{i+1}} (\xi_{i+1} - x)^k |f'(x)|^k\, dx.$$

Apply a Taylor expansion to $f'$ at the point $g_i$. Using (3.2), by means of the same arguments as above, we get

$$\left| \int_{f(t)}^{f(s)} \frac{|\tilde{U}_n(a) - g(a)|^k}{|g'(a)|^{k-1}}\, da \right.$$

$$\left. - \frac{1}{k+1} \sum_{i=0}^{m} |f'(g_i)|^k \{(\xi_{i+1} - g_i)^{k+1} - (\xi_{i+1} - g_{i+1})^{k+1}\} \right|$$

$$(3.4) \quad \leq C_1 \int_{g_i}^{g_{i+1}} (\xi_{i+1} - x)^k (x - g_i)\, dx$$

$$\leq C_1 (\xi_{i+1} - g_i) \int_{g_i}^{g_{i+1}} (\xi_{i+1} - x)^k\, dx$$

$$\leq \frac{C_1}{k+1} \{(\xi_{i+1} - g_i)^{k+2} - (\xi_{i+1} - g_{i+1})^{k+2}\}.$$

For the third integral in the statement of the lemma, similarly as before, again using (3.2), we can write

$$\int_{f(t)}^{f(s)} \frac{|\tilde{U}_n(a) - g(a)|^{k+1}}{|g'(a)|^k}\, da$$

$$(3.5) \quad = \sum_{i=0}^{m} \int_{g_i}^{g_{i+1}} |f'(g_i)|^{k+1} (\xi_{i+1} - x)^{k+1} \left| 1 + \frac{f''(\theta)}{f'(g_i)} (x - g_i) \right|^{k+1}$$

$$\geq \frac{C_2}{k+2} \sum_{i=0}^{m} \{(\xi_{i+1} - g_i)^{k+2} - (\xi_{i+1} - g_{i+1})^{k+2}\},$$

where $C_2 = (\inf |f'|/2)^{k+1}$.

Now let us define $\Delta$ as the difference between the first two integrals,

$$\Delta \stackrel{\text{def}}{=} \int_s^t |\tilde{f}_n(x) - f(x)|^k \, dx - \int_{f(t)}^{f(s)} \frac{|\tilde{U}_n(a) - g(a)|^k}{|g'(a)|^{k-1}} \, da.$$

By (3.3) and (3.4) and the fact that $\xi_0 = g_0$ and $\xi_{m+1} = g_{m+1}$, we find that

(3.6)
$$|\Delta| \leq D \sum_{i=0}^m (\xi_{i+1} - g_{i+1})^{k+1} ||f'(g_i)|^k - |f'(g_{i+1})|^k|$$
$$+ D \sum_{i=0}^m \{(\xi_{i+1} - g_i)^{k+2} - (\xi_{i+1} - g_{i+1})^{k+2}\},$$

where $D$ is a positive constant that depends only on the function $f$ and $k$. By a Taylor expansion, the first term on the right-hand side of (3.6) can be bounded by

$$D \sum_{i=0}^m (\xi_{i+1} - g_{i+1})^{k+1} |f'(g_i)|^k \left| 1 - \left| 1 + \frac{f''(\theta_i)(g_{i+1} - g_i)}{f'(g_i)} \right|^k \right|$$
$$\leq C_3 \sum_{i=0}^m (\xi_{i+1} - g_{i+1})^{k+1}(g_{i+1} - g_i)$$
$$\leq C_3 \sum_{i=0}^m \{(\xi_{i+1} - g_i)^{k+2} - (\xi_{i+1} - g_{i+1})^{k+2}\},$$

with $C_3$ depending only on $f$ and $k$, where we also use (3.2), the fact that $g_i < g_{i+1} < \xi_{i+1}$, and that according to (3.1), we have that $(g_{i+1} - g_i) \times \sup |f''|/\inf |f'| < \frac{1}{2}$. Putting things together and using (3.5), we find that

$$|\Delta| \leq C_4 \sum_{i=0}^m \{(\xi_{i+1} - g_i)^{k+2} - (\xi_{i+1} - g_{i+1})^{k+2}\}$$
$$\leq C_5 \int_{f(t)}^{f(s)} \frac{|\tilde{U}_n(a) - g(a)|^{k+1}}{|g'(a)|^k} \, da,$$

where $C_5$ depends only on $f$ and $k$. This proves the lemma for case 1. For case 2 the proof is similar. $\square$

**4. Asymptotic normality of the $L_k$-error of $\hat{f}_n$.** We will apply Lemma 3.1 to the following cut-off version of $\hat{f}_n$:

(4.1)
$$\tilde{f}_n(t) = \begin{cases} f(0), & \text{if } \hat{f}_n(x) \geq f(0), \\ \hat{f}_n(x), & \text{if } f(1) \leq \hat{f}_n(x) < f(0), \\ f(1), & \text{if } \hat{f}_n(x) < f(1). \end{cases}$$

The next lemma shows that $\tilde{f}_n$ satisfies condition (3.1) with probability tending to one.

LEMMA 4.1. *Define the event*

$$A_n = \left\{ \sup_{x \in [0,1]} |\tilde{f}_n(x) - f(x)| \leq \frac{\inf_{x \in [0,1]} |f'(x)|^2}{2 \sup_{t \in [0,1]} |f''(x)|} \right\}.$$

*Then* $P\{A_n^c\} \to 0$.

PROOF. It is sufficient to show that $\sup |\tilde{f}_n(x) - f(x)|$ tends to zero. For this we can follow the line of reasoning in Section 5.4 of [5]. Similar to their Lemma 5.9, we derive from our Lemma A.1 that, for each $a \in (f(1), f(0))$,

$$P\{|U_n(a) - g(a)| \geq n^{-1/3} \log n\} \leq C_1 \exp\{-C_2 (\log n)^3\}.$$

By monotonicity of $U_n$ and the conditions of $f$, this means that there exists a constant $C_3 > 0$ such that

$$P\left\{ \sup_{a \in (f(1), f(0))} |U_n(a) - g(a)| \geq C_3 n^{-1/3} \log n \right\} \leq C_1 \exp\{-\tfrac{1}{2} C_2 (\log n)^3\}.$$

This implies that the maximum distance between successive points of jump of $\hat{f}_n$ is of the order $\mathcal{O}(n^{-1/3} \log n)$. Since both $\tilde{f}_n$ and $f$ are monotone and bounded by $f(0)$, this also means that the maximum distance between $\tilde{f}_n$ and $f$ is of the order $\mathcal{O}(n^{-1/3} \log n)$. □

The difference between the $L_k$-errors for $\hat{f}_n$ and $\tilde{f}_n$ is bounded as

$$\begin{aligned}
(4.2) \quad & \left| \int_0^1 |\hat{f}_n(x) - f(x)|^k \, dx - \int_0^1 |\tilde{f}_n(x) - f(x)|^k \, dx \right| \\
& \leq \int_0^{U_n(f(0))} |\hat{f}_n(x) - f(x)|^k \, dx + \int_{U_n(f(1))}^1 |\hat{f}_n(x) - f(x)|^k \, dx.
\end{aligned}$$

The next lemma shows that the integrals on the right-hand side are of negligible order.

LEMMA 4.2. *Let* $U_n$ *be defined in* (1.1). *Then* $\int_0^{U_n(f(0))} |\hat{f}_n(x) - f(x)|^k \, dx = o_p(n^{-(2k+1)/6})$, *and* $\int_{U_n(f(1))}^1 |\hat{f}_n(x) - f(x)|^k \, dx = o_p(n^{-(2k+1)/6})$.

PROOF. Consider the first integral, which can be bounded by

$$\begin{aligned}
(4.3) \quad & 2^k \int_0^{U_n(f(0))} |\hat{f}_n(x) - f(0)|^k \, dx + 2^k \int_0^{U_n(f(0))} |f(x) - f(0)|^k \, dx \\
& \leq 2^k \int_0^{U_n(f(0))} |\hat{f}_n(x) - f(0)|^k \, dx + \frac{2^k}{k+1} \sup |f'|^k U_n(f(0))^{k+1}.
\end{aligned}$$

Define the event $B_n = \{U_n(f(0)) \leq n^{-1/3} \log n\}$. Then $U_n(f(0))^{k+1} \mathbb{1}_{B_n} = o_p(n^{-(2k+1)/6})$. Moreover, according to Theorem 2.1 in [4], it follows that $P\{B_n^c\} \to 0$. Since for any $\eta > 0$,

$$P(n^{(2k+1)/6} |U_n(f(0))|^{k+1} \mathbb{1}_{B_n^c} > \eta) \leq P\{B_n^c\} \to 0,$$

this implies that the second term in (4.3) is of the order $o_p(n^{-(2k+1)/6})$. The first term in (4.3) can be written as

$$
\begin{aligned}
& 2^k \left( \int_0^{U_n(f(0))} |\hat{f}_n(x) - f(0)|^k \, dx \right) \mathbb{1}_{B_n} \\
& + 2^k \left( \int_0^{U_n(f(0))} |\hat{f}_n(x) - f(0)|^k \, dx \right) \mathbb{1}_{B_n^c},
\end{aligned}
\tag{4.4}
$$

where the second integral is of the order $o_p(n^{-(2k+1)/6})$ by the same reasoning as before. To bound the first integral in (4.4), we will construct a suitable sequence $(a_i)_{i=1}^m$, such that the intervals $(0, n^{-a_1}]$ and $(n^{-a_i}, n^{-a_{i+1}}]$, for $i = 1, 2, \ldots, m-1$, cover the interval $(0, U_n(f(0))]$, and such that the integrals over these intervals can be bounded appropriately. First of all let

$$1 > a_1 > a_2 > \cdots > a_{m-1} \geq 1/3 > a_m, \tag{4.5}$$

and let $z_0 = 0$ and $z_i = n^{-a_i}$, $i = 1, \ldots, m$, so that $0 < z_1 < \cdots < z_{m-1} \leq n^{-1/3} < z_m$. On the event $B_n$, for $n$ sufficiently large, the intervals $(0, n^{-a_1}]$ and $(n^{-a_i}, n^{-a_{i+1}}]$ cover $(0, U_n(f(0))]$. Hence, when we denote $J_i = [z_i \wedge U_n(f(0)), z_{i+1} \wedge U_n(f(0))]$, the first integral in (4.4) can be bounded by

$$\sum_{i=0}^{m-1} \left( \int_{J_i} (\hat{f}_n(x) - f(0))^k \, dx \right) \mathbb{1}_{B_n} \leq \sum_{i=0}^{m-1} (z_{i+1} - z_i) |\hat{f}_n(z_i) - f(0)|^k,$$

using that $\hat{f}_n$ is decreasing and the fact that $J_i \subset (0, U_n(f(0))]$, so that $\hat{f}_n(z_i) - f(0) \geq \hat{f}_n(x) - f(0) \geq 0$, for $x \in J_i$. It remains to show that

$$\sum_{i=0}^{m-1} (z_{i+1} - z_i) |\hat{f}_n(z_i) - f(0)|^k = o_p(n^{-(2k+1)/6}). \tag{4.6}$$

From [13] we have that

$$\hat{f}_n(0) \to f(0) \sup_{1 \leq j < \infty} \frac{j}{\Gamma_j} \tag{4.7}$$

in distribution, where $\Gamma_j$ are partial sums of standard exponential random variables. Therefore,

$$z_1 |\hat{f}_n(0) - f(0)|^k = O_p(n^{-a_1}). \tag{4.8}$$

According to Theorem 3.1 in [9], for $1/3 \leq \alpha < 1$

$$(4.9) \qquad n^{(1-\alpha)/2}(\hat{f}_n(n^{-\alpha}) - f(n^{-\alpha})) \to Z$$

in distribution, where $Z$ is a nondegenerate random variable. Since for any $i = 1, \ldots, m-1$ we have that $1/3 \leq a_i < 1$, it follows that

$$|\hat{f}_n(z_i) - f(0)| \leq |\hat{f}_n(z_i) - f(z_i)| + \sup|f'|z_i$$
$$= \mathcal{O}_p(n^{-(1-a_i)/2}) + \mathcal{O}_p(n^{-a_i}) = \mathcal{O}_p(n^{-(1-a_i)/2}).$$

This implies that, for $i = 1, \ldots, m-1$,

$$(4.10) \qquad (z_{i+1} - z_i)|\hat{f}_n(z_i) - f(0)|^k = \mathcal{O}_p(n^{-a_{i+1}-k(1-a_i)/2}).$$

Therefore, if we can construct a sequence $(a_i)$ satisfying (4.5), as well as

$$(4.11) \qquad a_1 > \frac{2k+1}{6},$$

$$(4.12) \qquad a_{i+1} + \frac{k(1-a_i)}{2} > \frac{2k+1}{6} \qquad \text{for all } i = 1, \ldots, m-1,$$

then (4.6) follows from (4.8) and (4.10). One may take

$$a_1 = \frac{2k+7}{12},$$
$$a_{i+1} = \frac{k(a_i - 1)}{2} + \frac{2k+3}{8} \qquad \text{for } i = 1, \ldots, m-1.$$

Since $k < 2.5$, it immediately follows that (4.11) and (4.12) are satisfied. To show that (4.5) holds, first note that $1 > a_1 > 1/3$, because $k < 2.5$. It remains to show that the described sequence strictly decreases and reaches $1/3$ in finitely many steps. As long as $a_i > 1/3$, it follows that

$$a_i - a_{i+1} = \frac{2-k}{2}a_i + \frac{2k-3}{8}.$$

When $k = 2$, this equals $1/8$. For $1 \leq k < 2$, use $a_i > 1/3$, to find that $a_i - a_{i+1} > 1/24$, and for $2 \leq k < 2.5$, use $a_i \leq a_1 = (2k+1)/7$, to find that $a_i - a_{i+1} \geq (k+1)(2.5 - k)/12$. This means that the sequence $(a_i)$ also satisfies (4.5), which proves (4.6). This completes the proof of the first integral in the statement of the lemma. The proof for the second integral is similar. □

We are now able to prove our main result concerning the asymptotic normality of the $L_k$-error, for $1 \leq k < 2.5$.

PROOF OF THEOREM 1.1.  First consider the difference

$$(4.13) \qquad \left| \int_0^1 |\hat{f}_n(x) - f(x)|^k \, dx - \int_{f(1)}^{f(0)} \frac{|U_n(a) - g(a)|^k}{|g'(a)|^{k-1}} \, da \right|,$$

which can be bounded by

$$(4.14) \quad \left| \int_0^1 |\hat{f}_n(x) - f(x)|^k \, dx - \int_0^1 |\tilde{f}_n(x) - f(x)|^k \, dx \right| + R_n,$$

where

$$R_n = \left| \int_0^1 |\tilde{f}_n(x) - f(x)|^k \, dx - \int_{f(1)}^{f(0)} \frac{|U_n(a) - g(a)|^k}{|g'(t)|^{k-1}} \, da \right|.$$

Let $A_n$ be the event defined in Lemma 4.1, so that $P\{A_n^c\} \to 0$. As in the proof of Lemma 4.2, this means that $R_n \mathbb{1}_{A_n^c} = o_p(n^{-(2k+1)/6})$. Note that on the event $A_n$, the function $\tilde{f}_n$ satisfies the conditions of Lemma 3.1, and that for any $a \in [f(1), f(0)]$,

$$U_n(a) = \sup\{t \in [0,1] : \hat{f}_n(t) > a\} = \sup\{t \in [0,1] : \tilde{f}_n(t) > a\} = \tilde{U}_n(a).$$

Moreover, we can construct a partition $[0, s_1], (s_1, s_2], \ldots, (s_l, 1]$ of $[0, 1]$ in such a way that, on each element of the partition, $\tilde{f}_n$ satisfies either condition 1 or condition 2 of Lemma 3.1. This means that we can apply Lemma 3.1 to each element of the partition. Putting things together, it follows that $R_n \mathbb{1}_{A_n}$ is bounded from above by

$$C \int_{f(1)}^{f(0)} \frac{|U_n(a) - g(a)|^{k+1}}{|g'(a)|^k} \, da.$$

Corollary 2.1 implies that this integral is of the order $\mathcal{O}_p(n^{-(k+1)/3})$, so that $R_n \mathbb{1}_{A_n} = o_p(n^{-(2k+1)/6})$. Finally, the first difference in (4.14) can be bounded as in (4.2), which means that, according to Lemma 4.2, it is of the order $o_p(n^{-(2k+1)/6})$. Together with Corollary 2.1, this implies that

$$n^{1/6} \left( n^{k/3} \int_0^1 |\hat{f}_n(x) - f(x)|^k \, dx - \mu_k^k \right) \to N(0, \sigma^2),$$

where $\sigma^2$ is defined in Theorem 2.1. An application of the $\delta$-method then yields that

$$n^{1/6} \left( n^{1/3} \left( \int_0^1 |\hat{f}_n(x) - f(x)|^k \, dx \right)^{1/k} - \mu_k \right)$$

converges to a normal random variable with mean zero and variance

$$\left\{ \frac{1}{k} (\mu_k^k)^{1/k-1} \right\}^2 \sigma^2 = \frac{\sigma^2}{k^2 \mu_k^{2k-2}} = \sigma_k^2. \qquad \square$$

**5. Asymptotic normality of a modified $L_k$-error for large $k$.** For large $k$, the inconsistency of $\hat{f}_n$ at zero starts to dominate the behavior of the $L_k$-error. The following lemma indicates that, for $k > 2.5$, the result of Theorem 1.1 does not hold. For $k > 3$, the $L_k$-error tends to infinity, whereas for $2.5 < k \leq 3$, we are only able to prove that the variance of the integral near zero tends to infinity. In the latter case, it is in principle possible that the behavior of the process $\hat{f}_n - f$ on $[0, z_n]$ depends on the behavior of the process on $[z_n, 1]$ in such a way that the variance of the whole integral stabilizes, but this seems unlikely. The proof of this lemma is transferred to the Appendix.

LEMMA 5.1. *Let $z_n = 1/(2nf(0))$. Then we have the following:*

(i) *If $k > 3$, then $n^{k/3} E \int_0^1 |\hat{f}_n(x) - f(x)|^k \, dx \to \infty$.*
(ii) *If $k > 2.5$, then $\mathrm{var}(n^{(2k+1)/6} \int_0^{z_n} |\hat{f}_n(x) - f(x)|^k \, dx) \to \infty$.*

Although Lemma 5.1 indicates that, for $k > 2.5$, the result Theorem 1.1 will not hold for the usual $L_k$-error, a similar result can be derived for a modified version. For $k \geq 2.5$, we will consider a modified $L_k$-error of the form

$$(5.1) \qquad n^{1/6}\left\{ n^{1/3}\left( \int_{n^{-\varepsilon}}^{1-n^{-\varepsilon}} |\hat{f}_n(x) - f(x)|^k \, dx \right)^{1/k} - \mu_k \right\},$$

where $\mu_k$ is the constant defined in Theorem 1.1. In this way, for suitable choices of $\varepsilon$ we avoid a region where the Grenander estimator is inconsistent in such a way that we are still able to determine its global performance.

We first determine for what values of $\varepsilon$ we cannot expect asymptotic normality of (5.1). First of all, for $\varepsilon > 1$, similar to the proof of Lemma 5.1, it follows that

$$\mathrm{var}\left( n^{(2k+1)/6} \int_{n^{-\varepsilon}}^{z_n} |\hat{f}_n(x) - f(x)|^k \, dx \right) \to \infty.$$

For $\varepsilon < 1/6$, in view of Lemma 3.1 and the Brownian approximation discussed in Section 2, we have that the expectation of

$$n^{1/6}\left\{ n^{k/3} \int_{n^{-\varepsilon}}^{1-n^{-\varepsilon}} |\hat{f}_n(x) - f(x)|^k \, dx - \mu_k^k \right\}$$

will behave as the expectation of

$$n^{1/6}\left\{ \int_{f(n^{-\varepsilon})}^{f(1-n^{-\varepsilon})} \frac{n^{k/3}|U_n^W(a) - g(a)|^k}{|g'(a)|^{k-1}} \, da - \mu_k^k \right\},$$

which, according to Lemmas 2.1 and A.5, is of the order $\mathcal{O}(n^{1/6-\varepsilon})$. Hence, we also cannot expect asymptotic normality of (5.1) for $\varepsilon < 1/6$. Finally, for $(k-1)/(3k-6) < \varepsilon < 1$, a more tedious argument, in the same spirit as the proof of Lemma 5.1, yields that

$$\mathrm{var}\bigg(n^{(2k+1)/6} \int_{n^{-\varepsilon}}^{2n^{-\varepsilon}} |\hat{f}_n(x) - f(x)|^k \, dx\bigg) \to \infty.$$

Hence, in order to obtain a proper limit distribution for (5.1) for $k \geq 2.5$, we will choose $\varepsilon$ between $1/6$ and $(k-1)/(3k-6)$.

To prove a result analogous to Theorem 1.1, we define another cut-off version of the Grenander estimator,

$$f_n^\varepsilon(x) = \begin{cases} f(n^{-\varepsilon}), & \text{if } \hat{f}_n(x) \geq f(n^{-\varepsilon}), \\ \hat{f}_n(x), & \text{if } f(1 - n^{-\varepsilon}) \leq \hat{f}_n(x) < f(n^{-\varepsilon}), \\ f(1 - n^{-\varepsilon}), & \text{if } \hat{f}_n(x) < f(1 - n^{-\varepsilon}), \end{cases}$$

and its inverse function

(5.2) $$U_n^\varepsilon(a) = \sup\{x \in [n^{-\varepsilon}, 1 - n^{-\varepsilon}] : \hat{f}_n(x) \geq a\},$$

for $a \in [f(1 - n^{-\varepsilon}), f(n^{-\varepsilon})]$. The next lemma is the analogue of Lemma 4.1.

LEMMA 5.2. *Define the event*

$$A_n^\varepsilon = \bigg\{\sup_{x \in [0,1]} |f_n^\varepsilon(x) - f(x)| \leq \frac{\inf_{x \in [0,1]} |f'(x)|^2}{2 \sup_{t \in [0,1]} |f''(x)|}\bigg\}.$$

*Then $P\{A_n^\varepsilon\} \to 1$.*

PROOF. It suffices to show that $\sup_{x \in [0,1]} |f_n^\varepsilon(x) - f(x)| \to 0$. Using the definition of $f_n^\varepsilon$, we can bound

(5.3)
$$\sup_{x \in [0,1]} |f_n^\varepsilon(x) - f(x)|$$
$$\leq \sup_{x \in [0,1]} |f_n^\varepsilon(x) - \tilde{f}_n(x)| + \sup_{x \in [0,1]} |\tilde{f}_n(x) - f(x)|.$$

The first term on the right-hand side of (5.3) is smaller than $\sup |f'| n^{-\varepsilon}$, which, together with Lemma 4.1, implies that $\sup_{x \in [0,1]} |f_n^\varepsilon(x) - f(x)| = o_p(n^{-1/6})$. □

Similar to (4.2), the difference between the modified $L_k$-errors for $\hat{f}_n$ and $f_n^\varepsilon$ is bounded as

$$\bigg|\int_{n^{-\varepsilon}}^{1-n^{-\varepsilon}} |\hat{f}_n(x) - f(x)|^k \, dx - \int_{n^{-\varepsilon}}^{1-n^{-\varepsilon}} |f_n^\varepsilon(x) - f(x)|^k \, dx\bigg|$$

$$(5.4) \qquad \leq \int_{n^{-\varepsilon}}^{U_n^\varepsilon(f(n^{-\varepsilon}))} |\hat{f}_n(x) - f(x)|^k \, dx$$

$$+ \int_{U_n^\varepsilon(f(1-n^{-\varepsilon}))}^{1-n^{-\varepsilon}} |\hat{f}_n(x) - f(x)|^k \, dx.$$

The next lemma is the analogue of Lemma 4.2 and shows that both integrals on the right-hand side are of negligible order.

LEMMA 5.3. *For $k \geq 2.5$ and $1/6 < \varepsilon < (k-1)/(3k-6)$, let $U_n^\varepsilon$ be defined in (5.2). Then*

$$\int_{n^{-\varepsilon}}^{U_n^\varepsilon(f(n^{-\varepsilon}))} |\hat{f}_n(x) - f(x)|^k \, dx = o_p(n^{-(2k+1)/6})$$

*and*

$$\int_{U_n^\varepsilon(f(1-n^{-\varepsilon}))}^{1-n^{-\varepsilon}} |\hat{f}_n(x) - f(x)|^k \, dx = o_p(n^{-(2k+1)/6}).$$

PROOF. Consider the first integral. Then similar to (4.3), we have that

$$(5.5) \qquad 2^k \int_{n^{-\varepsilon}}^{U_n^\varepsilon(f(n^{-\varepsilon}))} |\hat{f}_n(x) - f(n^{-\varepsilon})|^k \, dx$$

$$+ 2^k \int_{n^{-\varepsilon}}^{U_n^\varepsilon(f(n^{-\varepsilon}))} |f(n^{-\varepsilon}) - f(x)|^k \, dx$$

$$\leq 2^k \int_{n^{-\varepsilon}}^{U_n^\varepsilon(f(n^{-\varepsilon}))} |\hat{f}_n(x) - f(n^{-\varepsilon})|^k \, dx$$

$$+ \frac{2^k}{k+1} \sup |f'|^k (U_n^\varepsilon(f(n^{-\varepsilon})) - n^{-\varepsilon})^{k+1}.$$

If we define the event $B_n^\varepsilon = \{U_n^\varepsilon(f(n^{-\varepsilon})) - n^{-\varepsilon} \leq n^{-1/3} \log n\}$, then by similar reasoning as in the proof of Lemma 4.2, it follows that $(U_n^\varepsilon(f(n^{-\varepsilon})) - n^{-\varepsilon})^{k+1} = o_p(n^{-(2k+1)/6})$. The first integral on the right-hand side of (5.5) can be written as

$$\left( \int_{n^{-\varepsilon}}^{U_n^\varepsilon(f(n^{-\varepsilon}))} |\hat{f}_n(x) - f(n^{-\varepsilon})|^k \, dx \right) \mathbb{1}_{B_n}$$

$$+ \left( \int_{n^{-\varepsilon}}^{U_n^\varepsilon(f(n^{-\varepsilon}))} |\hat{f}_n(x) - f(n^{-\varepsilon})|^k \, dx \right) \mathbb{1}_{B_n^c},$$

where the second term is of the order $o_p(n^{-(2k+1)/6})$ by the same reasoning as before. To bound

$$(5.6) \qquad \left( \int_{n^{-\varepsilon}}^{U_n^\varepsilon(f(n^{-\varepsilon}))} |\hat{f}_n(x) - f(n^{-\varepsilon})|^k \, dx \right) \mathbb{1}_{B_n},$$

we distinguish between two cases:

(i) $1/6 < \varepsilon \leq 1/3$,
(ii) $1/3 < \varepsilon < (k-1)/(3k-6)$.

In case (i), the integral (5.6) can be bounded by $|\hat{f}_n(n^{-\varepsilon}) - f(n^{-\varepsilon})|^k n^{-1/3} \log n$. According to Theorem 3.1 in [9], for $0 < \alpha < 1/3$,

$$(5.7) \qquad n^{1/3}(\hat{f}_n(n^{-\alpha}) - f(n^{-\alpha})) \to |4f(0)f'(0)|^{1/3} V(0)$$

in distribution, where $V(0)$ is defined in (1.2). It follows that $|\hat{f}_n(n^{-\varepsilon}) - f(n^{-\varepsilon})| = \mathcal{O}_p(n^{-1/3})$ and, therefore, (5.6) is of the order $o_p(n^{-(2k+1)/6})$.

In case (ii), similar to Lemma 4.2, we will construct a suitable sequence $(a_i)_{i=1}^m$, such that the intervals $(n^{-a_i}, n^{-a_{i+1}}]$, for $i = 1, 2, \ldots, m-1$, cover the interval $(n^{-\varepsilon}, U_n(f(n^{-\varepsilon}))]$, and such that the integrals over these intervals can be bounded appropriately. First of all let

$$(5.8) \qquad \varepsilon = a_1 > a_2 > \cdots > a_{m-1} \geq 1/3 > a_m,$$

and let $z_i = n^{-a_i}$, $i = 1, \ldots, m$, so that $0 < z_1 < \cdots < z_{m-1} \leq n^{-1/3} < z_m$. Then, similar to the proof of Lemma 4.2, we can bound (5.6) as

$$\left(\int_{n^{-\varepsilon}}^{U_n^\varepsilon(f(n^{-\varepsilon}))} |\hat{f}_n(x) - f(n^{-\varepsilon})|^k \, dx\right) \mathbb{1}_{B_n} \leq \sum_{i=1}^{m-1} (z_{i+1} - z_i)|\hat{f}_n(z_i) - f(n^{-\varepsilon})|^k.$$

Since $1/3 \leq a_i \leq \varepsilon < 1$ for $i = 1, \ldots, m-1$, we can apply (4.9) and conclude that each term is of the order $\mathcal{O}_p(n^{-a_{i+1} - k(1-a_i)/2})$. Therefore, it suffices to construct a sequence $(a_i)$ satisfying (5.8), as well as

$$(5.9) \qquad a_{i+1} + \frac{k(1-a_i)}{2} > \frac{2k+1}{6} \qquad \text{for all } i = 1, \ldots, m-1.$$

One may take

$$a_1 = \varepsilon,$$
$$a_{i+1} = \frac{k(a_i - 1)}{2} + \frac{2k+1}{6} + \frac{1}{8}\left(\frac{k-1}{3(k-2)} - \varepsilon\right) \qquad \text{for } i = 1, \ldots, m-1.$$

Then (5.9) is satisfied and it remains to show that the described sequence strictly decreases and reaches $1/3$ in finitely many steps. This follows from the fact that $a_i \leq \varepsilon$ and $k \geq 2.5$, since in that case

$$a_i - a_{i+1} = \frac{k-2}{2}\left(\frac{k-1}{3(k-2)} - a_i\right) - \frac{1}{8}\left(\frac{k-1}{3(k-2)} - \varepsilon\right)$$
$$\geq \frac{4k-9}{8}\left(\frac{k-1}{3(k-2)} - \varepsilon\right) > 0.$$

As in the proof of Lemma 4.2, the argument for the second integral is similar. Now take $B_n^\varepsilon = \{1 - n^{-\varepsilon} - U_n^\varepsilon(f(1 - n^{-\varepsilon})) \leq n^{-1/3} \log n\}$. The case $1/6 <$ $\varepsilon \leq 1/3$ can be treated in the same way as before. For the case $1/3 < \varepsilon < (k-1)/(3k-6)$, we can use the same sequence $(a_i)$ as above, but now define $z_i = 1 - n^{-a_i}$, $i = 1, \ldots, m$, so that $1 > z_1 > \cdots > z_{m-1} \geq 1 - n^{-1/3} > z_m$. Then we are left with considering

$$\left(\int_{U_n^\varepsilon(f(1-n^{-\varepsilon}))}^{1-n^{-\varepsilon}} |f(1 - n^{-\varepsilon}) - \hat{f}_n(x)|^k \, dx\right) \mathbb{1}_{B_n}$$

$$\leq \sum_{i=1}^{m-1} (z_i - z_{i+1}) |f(1 - n^{-\varepsilon}) - \hat{f}_n(z_i)|^k.$$

As before, each term in the sum is of the order $\mathcal{O}_p(n^{-a_{i+1} - k(1 - a_i)/2})$, for $i = 1, \ldots, m - 1$. The sequence chosen above satisfies (5.9) and (5.8), which implies that the sum above is of the order $o_p(n^{-(2k+1)/6})$. $\square$

Apart from (5.4), we also need to bound the difference between integrals for $U_n$ and its cut-off version $U_n^\varepsilon$:

(5.10)
$$\left| \int_{f(1)}^{f(0)} \frac{|U_n(a) - g(a)|^k}{|g'(a)|^{k-1}} \, da - \int_{f(1-n^{-\varepsilon})}^{f(n^{-\varepsilon})} \frac{|U_n^\varepsilon(a) - g(a)|^k}{|g'(a)|^{k-1}} \, da \right|$$
$$\leq \int_{\tilde{f}_n(n^{-\varepsilon})}^{f(0)} \frac{|U_n(a) - g(a)|^k}{|g'(a)|^{k-1}} \, da + \int_{f(1)}^{\tilde{f}_n(1-n^{-\varepsilon})} \frac{|U_n(a) - g(a)|^k}{|g'(a)|^{k-1}} \, da.$$

The next lemma shows that both integrals on the right-hand side are of negligible order.

LEMMA 5.4. *For $k \geq 2.5$, let $1/6 < \varepsilon < (k-1)/(3k-6)$. Furthermore, let $U_n$ be defined in (1.1) and let $\tilde{f}_n$ be defined in (4.1). Then*

$$\int_{\tilde{f}_n(n^{-\varepsilon})}^{f(0)} \frac{|U_n(a) - g(a)|^k}{|g'(a)|^{k-1}} \, da = o_p(n^{-(2k+1)/6})$$

*and*

$$\int_{f(1)}^{\tilde{f}_n(1-n^{-\varepsilon})} \frac{|U_n(a) - g(a)|^k}{|g'(a)|^{k-1}} \, da = o_p(n^{-(2k+1)/6}).$$

PROOF. Consider the first integral and define the event $A_n = \{f(0) - \tilde{f}_n(n^{-\varepsilon}) < n^{-1/6}/\log n\}$. For $1/6 < \varepsilon \leq 1/3$, according to (5.7) we have

$$f(0) - \tilde{f}_n(n^{-\varepsilon}) \leq |\hat{f}_n(n^{-\varepsilon}) - f(0)|$$
$$\leq |\hat{f}_n(n^{-\varepsilon}) - f(n^{-\varepsilon})| + \sup|f'| n^{-\varepsilon}$$
$$= \mathcal{O}_p(n^{-1/3}) + \mathcal{O}(n^{-\varepsilon})$$
$$= o_p(n^{-1/6}/\log n).$$

This means that, if $1/6 < \varepsilon \leq 1/3$, the probability $P\{A_n^c\} \to 0$. For $1/3 < \varepsilon < 1$,

$$P\{A_n^c\} \leq P\{f(0) - \tilde{f}_n(n^{-\varepsilon}) > 0\}$$
$$\leq P\{\hat{f}_n(n^{-\varepsilon}) - f(n^{-\varepsilon}) < n^{-\varepsilon} \sup|f'|\} \to 0,$$

since according to (4.9), $\hat{f}_n(n^{-\varepsilon}) - f(n^{-\varepsilon})$ is of the order $n^{-(1-\varepsilon)/2}$. Next write the first integral as

(5.11)
$$\left(\int_{\tilde{f}_n(n^{-\varepsilon})}^{f(0)} \frac{|U_n(a) - g(a)|^k}{|g'(a)|^{k-1}} \, da\right) \mathbb{1}_{A_n}$$
$$+ \left(\int_{\tilde{f}_n(n^{-\varepsilon})}^{f(0)} \frac{|U_n(a) - g(a)|^k}{|g'(a)|^{k-1}} \, da\right) \mathbb{1}_{A_n^c}.$$

Similar to the argument used in Lemma 4.2, the second integral in (5.11) is of the order $o_p(n^{-(2k+1)/6})$. The expectation of the first integral is bounded by

$$E \int_{f(0)-n^{-1/6}/\log n}^{f(0)} \frac{|U_n(a) - g(a)|^k}{|g'(a)|^{k-1}} \, da$$
$$\leq n^{-k/3} C_1 \int_{f(0)-n^{-1/6}/\log n}^{f(0)} E|V_n^E(a)|^k \, da$$
$$\leq C_2 n^{(2k+1)/6}/\log n,$$

using Lemma A.1. The Markov inequality implies that the first term in (5.11) is of the order $o_p(n^{-(2k+1)/6})$. For the second integral the proof is similar. □

THEOREM 5.1. *Suppose conditions* (A1)–(A3) *of Theorem* 1.1 *are satisfied. Then for $k \geq 2.5$ and for any $\varepsilon$ such that $1/6 < \varepsilon < (k-1)/(3k-6)$,*

$$n^{1/6}\left\{n^{1/3}\left(\int_{n^{-\varepsilon}}^{1-n^{-\varepsilon}} |\hat{f}_n(x) - f(x)|^k \, dx\right)^{1/k} - \mu_k\right\}$$

*converges in distribution to a normal random variable with zero mean and variance $\sigma_k^2$, where $\mu_k$ and $\sigma_k^2$ are defined in Theorem* 1.1.

PROOF. As in the proof of Theorem 1.1, it suffices to show that the difference

$$\left|\int_{n^{-\varepsilon}}^{1-n^{-\varepsilon}} |\hat{f}_n(x) - f(x)|^k \, dx - \int_{f(1)}^{f(0)} \frac{|U_n(a) - g(a)|^k}{|g'(a)|^{k-1}} \, da\right|$$

is of the order $o_p(n^{-(2k+1)/6})$. We can bound this difference by

$$(5.12) \quad \left| \int_{n^{-\varepsilon}}^{1-n^{-\varepsilon}} |\hat{f}_n(x) - f(x)|^k \, dx - \int_{n^{-\varepsilon}}^{1-n^{-\varepsilon}} |f_n^\varepsilon(x) - f(x)|^k \, dx \right|$$

$$(5.13) \quad + \left| \int_{f(1)}^{f(0)} \frac{|U_n(a) - g(a)|^k}{|g'(a)|^{k-1}} \, da - \int_{f(1-n^{-\varepsilon})}^{f(n^{-\varepsilon})} \frac{|U_n^\varepsilon(a) - g(a)|^k}{|g'(a)|^{k-1}} \, da \right|$$

$$(5.14) \quad + \left| \int_{n^{-\varepsilon}}^{1-n^{-\varepsilon}} |f_n^\varepsilon(x) - f(x)|^k \, dx - \int_{f(1-n^{-\varepsilon})}^{f(n^{-\varepsilon})} \frac{|U_n^\varepsilon(a) - g(a)|^k}{|g'(a)|^{k-1}} \, da \right|.$$

Differences (5.12) and (5.13) can be bounded as in (5.4) and (5.10), so that Lemmas 5.3 and 5.4 imply that these terms are of the order $o_p(n^{-(2k+1)/6})$. Finally, Lemma 3.1 implies that (5.14) is bounded by

$$\int_{f(1-n^{-\varepsilon})}^{f(n^{-\varepsilon})} \frac{|U_n^\varepsilon(a) - g(a)|^{k+1}}{|g'(a)|^k} \, da.$$

Write the integral as

$$\int_{f(1)}^{f(0)} \frac{|U_n(a) - g(a)|^{k+1}}{|g'(a)|^k} \, da$$

$$- \left( \int_{f(1)}^{f(0)} \frac{|U_n(a) - g(a)|^{k+1}}{|g'(a)|^k} \, da - \int_{f(1-n^{-\varepsilon})}^{f(n^{-\varepsilon})} \frac{|U_n^\varepsilon(a) - g(a)|^{k+1}}{|g'(a)|^k} \, da \right).$$

Then Corollary 2.1 and Lemma 5.4 imply that both terms are of the order $o_p(n^{-(2k+1)/6})$. This proves the theorem. □

## APPENDIX

The proofs in Section 2 follow the same line of reasoning as in [4]. Since we will frequently use results from this paper, we state them for easy reference. First, the tail probabilities of $V_n^J$ have a uniform exponential upper bound.

LEMMA A.1. *For $J = E, B, W$, let $V_n^J$ be defined by (2.2). Then there exist constants $C_1, C_2 > 0$ depending only on $f$, such that for all $n \geq 1$, $a \in (f(1), f(0))$ and $x > 0$, $P\{|V_n^J(a)| \geq x\} \leq C_1 \exp(-C_2 x^3)$.*

Properly normalized versions of $V_n^J(a)$ converge in distribution to $\xi(c)$ defined in (1.3). To be more precise, for $a \in (f(1), f(0))$ define $\phi_1(a) = |f'(g(a))|^{2/3}(4a)^{-1/3}$, $\phi_2(a) = (4a)^{1/3}|f'(g(a))|^{1/3}$ and

$$(A.1) \qquad V_{n,a}^J(c) = \phi_1(a) V_n^J(a - \phi_2(a) c n^{-1/3}),$$

for $J = E, B, W$. Then we have the following property.

LEMMA A.2. *For $J = E, B, W$, integer $d \geq 1$, $a \in (f(1), f(0))$ and $c \in J_n(a)^d$, we have joint distributional convergence of $(V_{n,a}^J(c_1), \ldots, V_{n,a}^J(c_d))$ to the random vector $(\xi(c_1), \ldots, \xi(c_d))$.*

Due to the fact that Brownian motion has independent increments, the process $V_n^W$ is mixing.

LEMMA A.3. *The process $\{V_n^W(a) : a \in (f(1), f(0))\}$ is strong mixing with mixing function $\alpha_n(d) = 12 e^{-C_3 n d^3}$, where the constant $C_3 > 0$ depends only on $f$.*

As a direct consequence of Lemma A.3 we have the following lemma, which is a slight extension of Lemma 4.1 in [4].

LEMMA A.4. *Let $l$ and $m$ be fixed such that $l + m > 0$ and let $h$ be a continuous function. Define*

$$c_h = 2 \int_0^1 (4f(x))^{(2l+2m+1)/3} |f'(x)|^{(4-4l-4m)/3} h(f(x))^2 \, dx.$$

*Then*

$$\operatorname{var}\left(n^{1/6} \int_{f(1)}^{f(0)} V_n^W(a)^l |V_n^W(a)|^m h(a) \, da\right)$$

$$\to c_h \int_0^\infty \operatorname{cov}(\xi(0)^l |\xi(0)|^m, \xi(c)^l |\xi(c)|^m) \, dc$$

*as $n \to \infty$.*

PROOF. The proof runs along the lines of the proof of Lemma 4.1 in [4]. We first have that

$$\operatorname{var}\left(n^{1/6} \int_{f(1)}^{f(0)} V_n^W(a)^l |V_n^W(a)|^m h(a) \, da\right)$$

$$= -2 \int_{f(1)}^{f(0)} \int_0^{n^{1/3} \phi_2(a)^{-1}(a-f(0))} (4a)^{(2l+2m+1)/3} |g'(a)|^{(4(l+m)-1)/3}$$

$$\times h(a) h(a - \phi_2(a) n^{-1/3} c)$$

$$\times \operatorname{cov}(V_{n,a}^W(0)^l |V_{n,a}^W(0)|^m, V_{n,a}^W(c)^l |V_{n,a}^W(c)|^m) \, dc \, da.$$

According to Lemma A.1, for $a$ and $c$ fixed, the sequence $V_{n,a}^W(c)^l |V_{n,a}^W(c)|^m$ is uniformly integrable. Hence, by Lemma A.2 the moments of $(V_{n,a}^W(0)^l |V_{n,a}^W(0)|^m, V_{n,a}^W(c)^l |V_{n,a}^W(c)|^m)$ converge to the corresponding moments of $(\xi(0)^l |\xi(0)|^m, \xi(c)^l |\xi(c)|^m)$. Again, Lemma A.1 and the fact that $l + m >$

0 yield that $E|V_{n,a}^W(0)|^{3(l+m)}$ and $E|V_{n,a}^W(c)|^{3(l+m)}$ are bounded uniformly in $n, a$ and $c$. Together with Lemma A.3 and Lemma 3.2 in [4], this yields that

$$|\text{cov}(V_{n,a}^W(0)^l|V_{n,a}^W(0)|^m, V_{n,a}^W(c)^l|V_{n,a}^W(c)|^m)| \leq D_1 e^{-D_2|c|^3},$$

where $D_1$ and $D_2$ do not depend on $n$, $a$ and $c$. The lemma now follows from dominated convergence and stationarity of the process $\xi$. □

PROOF OF THEOREM 2.1.  Write

$$W_n^k(a) = \frac{|V_n^W(a)|^k - E|V_n^W(a)|^k}{|g'(a)|^{k-1}},$$

and for $d = f(0) - f(1)$, define $L_n = dn^{-1/3}(\log n)^3$, $M_n = dn^{-1/3}\log n$ and $N_n = [d(L_n + M_n)^{-1}]$, where $[x]$ denotes the integer part of $x$. We divide the interval $(f(1), f(0))$ into $2N_n + 1$ blocks of alternating length,

$$A_j = (f(1) + (j-1)(L_n + M_n), f(1) + (j-1)(L_n + M_n) + L_n],$$
$$B_j = (f(1) + (j-1)(L_n + M_n) + L_n, f(1) + j(L_n + M_n)],$$

where $j = 1, \ldots, N_n$. Now write $T_{n,k} = S'_{n,k} + S''_{n,k} + R_{n,k}$, where

$$S'_{n,k} = n^{1/6} \sum_{j=1}^{N_n} \int_{A_j} W_n^k(a)\,da, \qquad S''_{n,k} = n^{1/6} \sum_{j=1}^{N_n} \int_{B_j} W_n^k(a)\,da,$$

$$R_{n,k} = n^{1/6} \int_{f(1)+N_n(L_n+M_n)}^{f(0)} W_n^k(a)\,da.$$

From here on the proof is completely the same as the proof of Theorem 4.1 in [4]. Therefore, we omit all specific details and only give a brief outline of the argument. Lemmas A.1 and A.3 imply that all moments of $W_n^k(a)$ are bounded uniformly in $a$ and that $E|W_n^k(a)W_n^k(b)| \leq D_1 \exp(-D_2 n|b-a|^3)$. This is used to ensure that $ER_n^2 \to 0$ and that the contribution of the small blocks is negligible: $E(S''_{n,k})^2 \to 0$. We then only have to consider the contribution over the big blocks. When we denote

$$Y_j = n^{1/6} \int_{A_j} W_n^k(a)\,da \quad \text{and} \quad \sigma_n^2 = \text{var}\left(\sum_{j=1}^{N_n} Y_j\right),$$

we find that

$$\left| E \exp\left\{ \frac{iu}{\sigma_n} \sum_{j=1}^{N_n} Y_j \right\} - \prod_{j=1}^{N_n} E \exp\left\{ \frac{iu}{\sigma_n} Y_j \right\} \right| \leq 4(N_n - 1) \exp(-C_3 n M_n^3) \to 0,$$

where $C_3 > 0$ depends only on $f$. This means that we can apply the central limit theorem to independent copies of $Y_j$. Since the moments of $|W_n^k(a)|$ are uniformly bounded, we have that, for each $\varepsilon > 0$,

$$\frac{1}{\sigma_n^2} \sum_{j=1}^{N_n} E Y_j^2 \mathbb{1}_{\{|Y_j| > \varepsilon \sigma_n\}} \leq \frac{1}{\varepsilon \sigma_n^3} N_n \sup_{1 \leq k \leq N_n} E|Y_j|^3 = \mathcal{O}(\sigma_n^{-3} n^{-1/6} (\log n)^6).$$

By similar computations as in the proof of Theorem 4.1 in [4], we find that $\sigma_n^2 = \text{var}(T_{n,k}) + \mathcal{O}(1)$, and application of Lemma A.4 yields that $\sigma_n^2 \to \sigma^2$. This implies that the $Y_j$'s satisfy the Lindeberg condition, which proves the theorem.

□

In order to prove Lemma 2.1, we first prove the following lemma.

LEMMA A.5. *Let $V_n^W$ be defined by (2.2) and let $V(0)$ be defined by (1.2). Then for $k \geq 1$, and for all $a$ such that*

(A.2) $$n^{1/3} \{F(g(a)) \wedge (1 - F(g(a)))\} \geq \log n,$$

*we have*

$$E|V_n^W(a)|^k = E|V(0)|^k \frac{(4a)^{k/3}}{|f'(g(a))|^{2k/3}} + \mathcal{O}(n^{-1/3} (\log n)^{k+3}),$$

*where the term $\mathcal{O}(n^{-1/3} (\log n)^{k+3})$ is uniform in all $a$ satisfying (A.2).*

PROOF. The proof relies on the proof of Corollary 3.2 in [4]. There it is shown that, if we define $H_n(y) = n^{1/3} \{H(F(g(a)) + n^{-1/3} y) - g(a)\}$, with $H$ being the inverse of $F$, and

$$V_{n,b} = \sup\{y \in [-n^{1/3} F(g(a)), n^{1/3} (1 - F(g(a)))] : W(y) - by^2 \text{ is maximal}\},$$

with $b = |f'(g(a))|/(2a^2)$, then for the event $A_n = \{|V_n^W(a)| \leq \log n, |H_n(V_{n,b})| \leq \log n\}$, one has that $P\{A_n^c\}$ is of the order $\mathcal{O}(e^{-C(\log n)^3})$, which then implies that

$$\sup_{a \in (f(1), f(0))} E|V_n^W(a) - H_n(V_{n,b})| = \mathcal{O}(n^{-1/3} (\log n)^4).$$

Similarly, together with an application of the mean value theorem, this yields

(A.3) $$\sup_{a \in (f(1), f(0))} E||V_n^W(a)|^k - |H_n(V_{n,b})|^k| = \mathcal{O}(n^{-1/3} (\log n)^{3+k}).$$

Note that, by definition, the $\arg\max V_{n,b}$ closely resembles the $\arg\max V_b(0)$, where

(A.4) $$V_b(c) = \arg\max_{t \in \mathbb{R}} \{W(t) - b(t-c)^2\}.$$

Therefore, we write

(A.5) $\quad E|H_n(V_{n,b})|^k = E|H_n(V_b(0))|^k + E(|H_n(V_{n,b})|^k - |H_n(V_b(0))|^k).$

Since by Brownian scaling $V_b(c)$ has the same distribution as $b^{-2/3}V(cb^{2/3})$, where $V$ is defined in (1.2), together with the conditions on $f$, we find that

$$E|H_n(V_b(0))|^k = a^{-k}E|V_b(0)|^k + \mathcal{O}(n^{-1/3})$$
$$= \frac{(4a)^{k/3}}{|f'(g(a))|^{2k/3}}E|V(0)|^k + \mathcal{O}(n^{-1/3}).$$

As in the proof of Corollary 3.2 in [4], $V_{n,b}$ can only be different from $V_b(0)$ with probability of the order $\mathcal{O}(e^{-(2/3)(\log n)^3})$. Hence, from (A.5), we conclude that

$$E|H_n(V_{n,b})|^k = \frac{(4a)^{k/3}}{|f'(g(a))|^{2k/3}}E|V(0)|^k + \mathcal{O}(n^{-1/3}).$$

Together with (A.3), this proves the lemma. □

PROOF OF LEMMA 2.1. The result immediately follows from Lemma A.5. The values of $a$ for which condition (A.2) does not hold give a contribution of the order $\mathcal{O}(n^{-1/3}\log n)$ to the integral $\int E|V_n^W(a)|^k\,da$, and finally,

$$\int_{f(1)}^{f(0)} \frac{(4a)^{k/3}}{|f'(g(a))|^{2k/3}|g'(a)|^{k-1}}\,da = \int_0^1 (4f(x))^{k/3}|f'(x)|^{k/3}\,dx. \quad \square$$

PROOF OF LEMMA 2.2. The proof of the first statement relies on the proof of Corollary 3.3 in [4]. Here it is shown that, if for $a$ belonging to the set $J_n = \{a : \text{both } a \text{ and } a(1 - \xi_n n^{-1/2}) \in (f(1), f(0))\}$ we define

$$V_n^B(a, \xi_n) = V_n^B(a(1 - n^{-1/2}\xi_n)) + n^{1/3}\{g(a(1 - n^{-1/2}\xi_n)) - g(a)\},$$

then for the event $A_n = \{|\xi_n| \leq n^{1/6}, |V_n^W(a)| \leq \log n, |V_n^B(a,\xi_n)| \leq \log n\}$, one has that $P\{A_n^c\}$ is of the order $\mathcal{O}(e^{-C(\log n)^3})$. This implies that

$$\int_{a \in J_n} E|V_n^B(a, \xi_n) - V_n^W(a)|\,da = \mathcal{O}(n^{-1/3}(\log n)^3).$$

Hence, by using the same method as in proof of Lemma A.5, we obtain

$$\int_{a \in J_n} E||V_n^B(a, \xi_n)|^k - |V_n^W(a)|^k|\,da = \mathcal{O}(n^{-1/3}(\log n)^{k+2}).$$

From Lemma A.1, it also follows that $E|V_n^B(a)|^k = \mathcal{O}(1)$ and $E|V_n^W(a)|^k = \mathcal{O}(1)$, uniformly with respect to $n$ and $a \in (f(1), f(0))$. Hence, the contribution of the integrals over $[f(1), f(0)] \setminus J_n$ is negligible, and it remains to show that

(A.6) $\quad n^{1/6} \int_{a \in J_n} \{|V_n^B(a, \xi_n)|^k - |V_n^B(a)|^k\}\,da = o_p(1).$

For $k = 1$ this is shown in the proof of Corollary 3.3 in [4], so we may assume that $k > 1$. Completely similar to the proof in the case $k = 1$, we first obtain

$$n^{1/6} \int_{a \in J_n} \{|V_n^B(a, \xi_n)|^k - |V_n^B(a)|^k\} \, da$$
$$= n^{1/6} \int_{f(1)}^{f(0)} \{|V_n^B(a) - ag'(a)\xi_n n^{-1/6}|^k - |V_n^B(a)|^k\} \, da + \mathcal{O}_p(n^{-1/3}).$$

Let $\varepsilon > 0$ and write $\Delta_n(a) = ag'(a)\xi_n n^{-1/6}$. Then the first term on the right-hand side equals

(A.7) $\quad n^{1/6} \int_{f(1)}^{f(0)} \{|V_n^B(a) - \Delta_n(a)|^k - |V_n^B(a)|^k\} \mathbb{1}_{[0,\varepsilon]}(|V_n^B(a)|) \, da$

(A.8) $\quad + n^{1/6} \int_{f(1)}^{f(0)} \{|V_n^B(a) - \Delta_n(a)|^k - |V_n^B(a)|^k\} \mathbb{1}_{(\varepsilon,\infty)}(|V_n^B(a)|) \, da.$

First consider the term (A.7). When $|V_n^B(a)| < 2|\Delta_n(a)|$, we can write

$$||V_n^B(a) - \Delta_n(a)|^k - |V_n^B(a)|^k| \leq 3^k |\Delta_n(a)|^k + 2^k |\Delta_n(a)|^k$$
$$\leq (3^k + 2^k)|ag'(a)\xi_n|^k n^{-k/6}.$$

When $|V_n^B(a)| \geq 2|\Delta_n(a)|$, we have

$$||V_n^B(a) - \Delta_n(a)|^k - |V_n^B(a)|^k| = k|\theta|^{k-1}|ag'(a)\xi_n|n^{-1/6},$$

where $\theta$ is between $|V_n^B(a)| \leq \varepsilon$ and $|V_n^B(a) - \Delta_n(a)| \leq \frac{3}{2}\varepsilon$. Using that $\xi_n$ and $V_n^B$ are independent, the expectation of (A.7) is bounded from above by

$$C_1 \varepsilon^{k-1} E|\xi_n| \int_{f(1)}^{f(0)} |ag'(a)| P\{|V_n^B(a)| \leq \varepsilon\} \, da + \mathcal{O}_p(n^{-(k-1)/6}),$$

where $C_1 > 0$ depends only on $f$ and $k$. Hence, since $k > 1$, we find that

(A.9)
$$\limsup_{n \to \infty} n^{1/6} \int_{f(1)}^{f(0)} \{|V_n^B(a) - ag'(a)\xi_n n^{-1/6}|^k - |V_n^B(a)|^k\}$$
$$\times \mathbb{1}_{[0,\varepsilon]}(|V_n^B(a)|) \, da$$

is bounded from above by $C_2 \varepsilon^{k-1}$, where $C_2 > 0$ depends only on $f$ and $k$. Letting $\varepsilon \downarrow 0$ and using that $k > 1$ then yields that (A.7) tends to zero.

The term (A.8) is equal to

(A.10)
$$\int_{f(1)}^{f(0)} \frac{-2\xi_n ag'(a) V_n^B(a) + (ag'(a)\xi_n)^2 n^{-1/6}}{|V_n^B(a) - \Delta_n(a)| + |V_n^B(a)|}$$
$$\times k\theta(a)^{k-1} \mathbb{1}_{(\varepsilon,\infty)}(|V_n^B(a)|) \, da,$$

where $\theta(a)$ is between $|V_n^B(a) - \Delta_n(a)|$ and $|V_n^B(a)|$. Note that for $|V_n^B(a)| > \varepsilon$,

$$\left| \frac{2V_n^B(a)}{|V_n^B(a) - \Delta_n(a)| + |V_n^B(a)|} - \frac{V_n^B(a)}{|V_n^B(a)|} \right| \leq \frac{|ag'(a)n^{-1/6}\xi_n|}{\varepsilon} = \mathcal{O}_p(n^{-1/6})$$

uniformly in $a \in (f(1), f(0))$, so that (A.10) is equal to

$$-k\xi_n \int_{f(1)}^{f(0)} ag'(a) V_n^B(a) |V_n^B(a)|^{k-2} \mathbb{1}_{(\varepsilon,\infty)}(|V_n^B(a)|) \, da$$
$$+ k\xi_n \int_{f(1)}^{f(0)} ag'(a) \frac{V_n^B(a)}{|V_n^B(a)|} (|V_n^B(a)|^{k-1} - \theta(a)^{k-1}) \mathbb{1}_{(\varepsilon,\infty)}(|V_n^B(a)|) \, da$$
$$+ \mathcal{O}_p(n^{-1/6}).$$

We have that

$$||V_n^B(a)|^{k-1} - \theta(a)^{k-1}| \leq |V_n^B(a)|^{k-1} \left| \left|1 - \frac{\Delta_n(a)}{V_n^B(a)}\right|^{k-1} - 1 \right| = \mathcal{O}_p(n^{-1/6}),$$

where the $\mathcal{O}$-term is uniform in $a$. This means that (A.10) is equal to

(A.11) $\quad -k\xi_n \int_{f(1)}^{f(0)} ag'(a) V_n^B(a) |V_n^B(a)|^{k-2} \, da$

$\qquad\qquad + k\xi_n \int_{f(1)}^{f(0)} ag'(a) \operatorname{sign}(V_n^B(a)) |V_n^B(a)|^{k-1} \mathbb{1}_{[0,\varepsilon)}(|V_n^B(a)|) \, da$

(A.12) $\qquad\qquad + \mathcal{O}_p(n^{-1/6}).$

The integral in (A.12) is of the order $\mathcal{O}(\varepsilon^{k-1})$, whereas $E\xi_n^2 = 1$. Since $k > 1$, this means that, after letting $\varepsilon \downarrow 0$, (A.12) tends to zero. For (A.11), let $S_n^B(a) = ag'(a) V_n^B(a) |V_n^B(a)|^{k-2}$ and consider $E(\int S_n^B(a) \, da)^2 = \operatorname{var}(\int S_n^B(a) \, da) + (E \int S_n^B(a) \, da)^2$. Then, since according to Lemma A.1 all moments of $|S_n^B(a)|$ are bounded uniformly in $a$, we find by dominated convergence and Lemma A.2 that

$$\lim_{n \to \infty} E \int S_n^B(a) \, da = \int \frac{a|g'(a)|}{(\phi_1(a))^k} (E\xi(0)|\xi(0)|^{k-2}) \, da = 0,$$

because the distribution of $\xi(0)$ is symmetric. Applying Lemma A.4 with $l = 1$, $m = k - 2$ and $h(a) = ag'(a)$, we obtain $\operatorname{var}(\int S_n^B(a) \, da) = O(n^{-1/3})$. We conclude that (A.10) tends to zero in probability. This proves the first statement of the lemma.

The proof of the second statement relies on the proof of Corollary 3.1 in [4]. There it is shown that, for the event $A_n = \{|V_n^B(a)| < \log n, |V_n^E(a)| <$ $\log n\}$ one has that $P\{A_n^c\}$ is of the order $\mathcal{O}(e^{-C(\log n)^3})$. Furthermore, if $K_n = \{\sup_t |E_n(t) - B_n(F(t))| \leq n^{-1/2}(\log n)^2\}$, then $P\{K_n\} \to 1$ and

$$\text{(A.13)} \qquad E||V_n^E(a)| - |V_n^B(a)||\mathbb{1}_{A_n \cap K_n} = \mathcal{O}(n^{-1/3}(\log n)^3)$$

uniformly in $a \in (f(1), f(0))$. By the mean value theorem, together with (A.13), we now have that

$$E||V_n^E(a)|^k - |V_n^B(a)|^k|\mathbb{1}_{K_n}$$
$$\leq k(\log n)^{k-1} E||V_n^E(a)| - |V_n^B(a)||\mathbb{1}_{A_n \cap K_n} + 2n^{k/3} P\{A_n^c\}$$
$$= \mathcal{O}(n^{-1/3}(\log n)^{k+2}) + \mathcal{O}(n^{k/3} e^{-C(\log n)^3}).$$

This proves the lemma. □

This completes the proofs needed in Section 2 to obtain a central limit theorem for the scaled $L_k$-distance between $U_n$ and $g$ (Corollary 2.1). The remainder of this appendix is devoted to the proof of Lemma 5.1, which indicates that a central limit theorem for the $L_k$-distance between $\hat{f}_n$ and $f$ is not possible when $k > 2.5$. For this we need the following lemma.

LEMMA A.6. *Let $k \geq 2.5$ and $z_n = 1/(2nf(0))$. Then there exist $0 < a_1 < b_1 < a_2 < b_2 < \infty$, such that, for $i = 1, 2$,*

$$\liminf_{n \to \infty} P\left\{n \int_0^{z_n} |\hat{f}_n(x) - f(x)|^k \, dx \in [a_i, b_i]\right\} > 0.$$

PROOF. Consider the event $A_n = \{X_i \geq z_n, \text{ for all } i = 1, 2, \ldots, n\}$. Then it follows that $P\{A_n\} \to 1/\sqrt{e} > 1/2$. Since on the event $A_n$ the estimator $\hat{f}_n$ is constant on the interval $[0, z_n]$, for any $a_i > 0$ we have

$$P\left\{n \int_0^{z_n} |\hat{f}_n(x) - f(x)|^k \, dx \in [a_i, b_i]\right\}$$
$$\geq P\left\{\left(n \int_0^{z_n} |\hat{f}_n(0) - f(x)|^k \, dx\right) \mathbb{1}_{A_n} \in [a_i, b_i]\right\}$$
$$= P\left\{\left(\frac{|\hat{f}_n(0) - f(0)|^k}{2f(0)} + R_n\right) \mathbb{1}_{A_n} \in [a_i, b_i]\right\},$$

where

$$R_n = n \int_0^{z_n} k\theta_n(x)^{k-1}(|\hat{f}_n(0) - f(x)| - |\hat{f}_n(0) - f(0)|) \, dx,$$

with $\theta_n(x)$ between $|\hat{f}_n(0) - f(x)|$ and $|\hat{f}_n(0) - f(0)|$. Using (4.7), we obtain that $R_n$ is of the order $\mathcal{O}_p(n^{-1})$ and, therefore,

$$\frac{|\hat{f}_n(0) - f(0)|^k}{2f(0)} + R_n \to \frac{f(0)^{k-1}}{2}\left|\sup_{1 \leq j < \infty} \frac{j}{\Gamma_j} - 1\right|^k$$

in distribution. Now choose $0 < a_1 < b_1 < a_2 < b_2 < \infty$ such that, for $i = 1, 2$,

$$P\left\{\frac{f(0)^{k-1}}{2}\left|\sup_{1\le j<\infty}\frac{j}{\Gamma_j}-1\right|^k \in [a_i, b_i]\right\} > 1 - 1/\sqrt{e}.$$

Then for $i = 1, 2$ we find

$$P\left\{n\int_0^{z_n} |\hat{f}_n(x) - f(x)|^k\, dx \in [a_i, b_i]\right\}$$
$$\ge P\left\{\left(\frac{|\hat{f}_n(0) - f(0)|^k}{2f(0)} + R_n\right) \in [a_i, b_i]\right\} - P\{A_n^c\},$$

which converges to a positive value. $\square$

PROOF OF LEMMA 5.1. Take $0 < a_1 < b_1 < a_2 < b_2 < \infty$ as in Lemma A.6, and let $A_{ni}$ be the event

$$A_{ni} = \left\{n\int_0^{z_n} |\hat{f}_n(x) - f(x)|^k\, dx \in [a_i, b_i]\right\}.$$

Then

$$n^{k/3} E \int_0^1 |\hat{f}_n(x) - f(x)|^k\, dx \ge n^{k/3} E \int_0^{z_n} |\hat{f}_n(x) - f(x)|^k\, dx\, \mathbb{1}_{A_{n1}}$$
$$\ge a_1 n^{(k-3)/3} P\{A_{n1}\}.$$

Since according to Lemma A.6 $P\{A_{n1}\}$ tends to a positive constant, this proves (i).

For (ii), write $X_n = n\int_0^{z_n} |\hat{f}_n(x) - f(x)|^k\, dx$, and define $B_n = \{EX_n \ge (a_2 + b_1)/2\}$. Then

$$\operatorname{var}(X_n) \ge E(X_n - EX_n)^2 \mathbb{1}_{A_{n1}\cap B_n} + E(X_n - EX_n)^2 \mathbb{1}_{A_{n2}\cap B_n^c}$$
$$\ge \tfrac{1}{4}(a_2 - b_1)^2 P\{A_{n1}\}\mathbb{1}_{B_n} + \tfrac{1}{4}(a_2 - b_1)^2 P\{A_{n2}\}\mathbb{1}_{B_n^c}$$
$$\ge \tfrac{1}{4}(a_2 - b_1)^2 \min(P\{A_{n1}\}, P\{A_{n2}\}).$$

Hence, according to Lemma A.6,

$$\liminf_{n\to\infty} \operatorname{var}\left(n^{(2k+1)/6}\int_0^{z_n} |\hat{f}_n(x) - f(x)|^k\, dx\right)$$
$$\ge \liminf_{n\to\infty} n^{(2k-5)/3}\tfrac{1}{4}(a_2 - b_1)^2 \min(P\{A_{n1}\}, P\{A_{n2}\}) = \infty. \quad \square$$

EURANDOM
P.O. BOX 513-5600 MB
EINDHOVEN
THE NETHERLANDS
E-MAIL: kulikov@eurandom.nl

DEPARTMENT CROSS
FACULTY OF ELECTRICAL ENGINEERING,
  MATHEMATICS AND COMPUTER SCIENCE
DELFT UNIVERSITY OF TECHNOLOGY
MEKELWEG 4
2628 CD DELFT
THE NETHERLANDS
E-MAIL: h.p.lopuhaa@its.tudelft.nl